\documentclass[a4paper,11pt]{amsart}
\usepackage{amssymb,amsmath}
\usepackage{amsopn}
\usepackage{mathrsfs}
\usepackage{color}
\usepackage{wasysym}
\usepackage{vmargin}
\usepackage{vmargin}
\usepackage{color}
\usepackage{amscd}
\usepackage{verbatim}
\newtheorem{theorem}{Theorem}[section]

\numberwithin{equation}{section}

\newcommand\RR{{{\mathbb R}}}

\newcommand\SSS{{\mathbb S}}
\newcommand{\rr}{\mathbb{R}}
\newcommand{\eps}{\varepsilon}
\newcommand{\nn}{\mathbb{N}}

\def\p{\partial}

\def\Id{\operatorname{Id}}

\def\N{\mathbb N}

\def\R{\mathbb R}

\def\val#1{\vert#1\vert}

\def\l2{L^2(\R^{n})}
\def\L2{L^2(\R^{2n})}

\def\vs{\vskip.3cm}

\def\mat22#1#2#3#4{\begin{pmatrix}#1&#2\\ #3&#4\end{pmatrix}}

\begin{document}

\title[On the ultra-analytic smoothing properties of the Landau equation]
{A remark on the ultra-analytic smoothing properties of the spatially homogeneous Landau equation}
\author{Y. Morimoto, K. Pravda-Starov  \& C.-J. Xu}
\date{\today}
\address{\noindent \textsc{Y. Morimoto, Graduate School of Human and Environmental Studies,
Kyoto University, Kyoto 606-8501, Japan}}
\email{morimoto@math.h.kyoto-u.ac.jp }
\address{\noindent \textsc{K. Pravda-Starov,
Universit\'e de Cergy-Pontoise,
D\'epartement de Math\'ematiques, CNRS UMR 8088,
95000 Cergy-Pontoise, France}}
\email{karel.pravda-starov@u-cergy.fr}
\address{\noindent \textsc{C.-J. Xu, School of Mathematics, Wuhan university 430072, Wuhan, P.R. China\\
 and  \\
 Universit\'e de Rouen, CNRS UMR 6085, D\'epartement de Math\'ematiques, 76801 Saint-Etienne du Rouvray, France}}
\email{Chao-Jiang.Xu@univ-rouen.fr}

\begin{abstract}
We consider the non-linear spatially homogeneous Landau equation with Maxwellian molecules in a close-to-equilibrium framework and show that the Cauchy problem for the fluctuation around the Maxwellian equilibrium distribution enjoys a Gelfand-Shilov regularizing effect in the class $S_{1/2}^{1/2}(\rr^d)$, implying the ultra-analyticity of both the fluctuation and its Fourier transform, for any positive time.
\end{abstract}

\maketitle

\section{Introduction}
In the work~\cite{LMPX3}, we consider the spatially homogeneous non-cutoff Boltzmann equation with Maxwellian molecules in a close-to-equilibrium framework and
study the smoothing properties of the Cauchy problem for the fluctuation around the Maxwellian equilibrium distribution. 
The Boltzmann equation describes the behavior of a dilute gas when the only interactions taken into account are binary collisions~\cite{17}. In the spatially homogeneous case with Maxwellian molecules, it reads as the equation
\begin{equation}\label{e1vv}
\begin{cases}
\partial_tf=Q(f,f),\\
f|_{t=0}=f_0,
\end{cases}
\end{equation}
for the density distribution of the particles $f=f(t,v) \geq 0$, $t \geq 0$, $v \in \rr^d$, with $d \geq 2$, where the non-linear term 
\begin{equation}\label{eq1vv}
Q(f, f)=\int_{\rr^d}\int_{\SSS^{d-1}}b\Big(\frac{v-v_*}{|v-v_*|} \cdot \sigma\Big)(f'_* f'-f_*f)d\sigma dv_*,
\end{equation}
stands for the Boltzmann collision operator whose cross section is a non-negative function satisfying to the assumption
\begin{equation}\label{sa1vv}
(\sin \theta)^{d-2}b(\cos \theta)  \substack{\\ \\ \approx \\ \theta \to 0_{+} }  \theta^{-1-2s},
\end{equation}
for some  $0 < s <1$. The notation $a\approx b$ means that $a/b$ is bounded from above and below by fixed positive constants. The term (\ref{sa1vv}) is not integrable in zero
$$\int_0^{\frac{\pi}{2}}(\sin \theta)^{d-2}b(\cos \theta)d\theta=+\infty.$$
This non-integrability plays a major role regarding the qualitative behaviour of the solutions to the Boltzmann equation and this feature is essential for the smoothing effect to be present, see the discussion in~\cite{LMPX3} and all the references herein.

In~\cite{LMPX3}, we consider the spatially homogeneous non-cutoff Boltzmann equation with Maxwellian molecules (\ref{e1vv}) in the radially symmetric case with initial density distributions 
\begin{equation}\label{kk1}
f_0=\mu_d+\sqrt{\mu_d}g_0, \quad g_0 \in L^2(\rr^d) \textrm{ radial},\quad  \|g_0\|_{L^2} \ll 1,
\end{equation}
close to the Maxwellian equilibrium distribution 
\begin{equation}\label{maxwe}
\mu_d(v)=(2\pi)^{-\frac{d}{2}}e^{-\frac{|v|^2}{2}}, \quad v \in \rr^d, \qquad Q(\mu_d,\mu_d)=0,
\end{equation}
where $|\cdot|$ is the Euclidean norm on $\rr^d$, in the physical 3-dimensional case $d=3$. The main result in~\cite{LMPX3} shows that the Cauchy problem for the fluctuation 
$$f=\mu_3+\sqrt{\mu_3}g,$$ 
around the Maxwellian equilibrium distribution 
\begin{equation}\label{bs1}
\begin{cases}
\partial_tg+\mathscr{L}g=\mu_3^{-1/2}Q(\sqrt{\mu_3}g,\sqrt{\mu_3}g),\\
g|_{t=0}=g_0 \in L^2(\rr^3),
\end{cases}
\end{equation}
where
$$\mathscr{L}g=-\mu_3^{-1/2}Q(\mu_3,\mu_3^{1/2}g)-\mu_3^{-1/2}Q(\mu_3^{1/2}g,\mu_3),$$
enjoys the same Gelfand-Shilov regularizing effect as the Cauchy problem defined by the evolution equation associated to the fractional harmonic oscillator
\begin{equation}\label{gel2vv}
\begin{cases}
\partial_tg+\mathcal{H}^sg=0,\\
g|_{t=0}=g_0 \in L^2(\rr^3),
\end{cases}\qquad \mathcal{H}=-\Delta_v+\frac{|v|^2}{4},
\end{equation}
where $0<s<1$ is the positive parameter appearing in the assumption (\ref{sa1vv}). More specifically, we prove that under the assumption (\ref{kk1}), the Cauchy problem (\ref{bs1}) admits a unique global radial solution $g \in L^{\infty}(\rr_t^+,L^2(\rr_v^3))$, which belongs to the Gelfand-Shilov class $S_{1/2s}^{1/2s}(\rr^3)$ for any positive time
\begin{equation}\label{rrr1}
\forall t >0, \quad g(t) \in S_{1/2s}^{1/2s}(\rr^3).
\end{equation}
The definition of the Gelfand-Shilov regularity is recalled in appendix (Section~\ref{regularity}).

In the present work, we study the spatially homogeneous Landau equation with Maxwellian molecules
\begin{equation}\label{e1bisggg}
\begin{cases}
\partial_tf=Q_L(f,f),\\
f|_{t=0}=f_0.
\end{cases}
\end{equation}
The Landau collision operator $Q_L(f,f)$ is understood as the limiting Boltzmann operator in the grazing collision limit asymptotic \cite{mou5,mou1,mou2,mou3,villani1}, when $s$ tends to 1 in the singularity assumption (\ref{sa1vv}).
In the physical 3-dimensional case, the linearized non-cutoff Boltzmann operator with Maxwellian molecules was actually showed to be equal to the fractional linearized Landau operator with Maxwellian molecules \cite{LMPX2} (Theorem~2.3),
$$\mathscr{L}=a(\mathcal{H},\Delta_{\SSS^2}) \mathscr{L}_L^s,$$
up to a positive bounded isomorphism on $L^2(\rr^3)$,
$$\exists c>0, \forall f \in L^2(\rr^3), \quad c\|f\|_{L^2}^2 \leq (a(\mathcal{H},\Delta_{\SSS^2})f,f)_{L^2} \leq \frac{1}{c}\|f\|_{L^2}^2,$$
commuting with the harmonic oscillator 
$\mathcal H=-\Delta_v+\frac{|v|^2}{4}$ 
and the Laplace-Beltrami operator
$$\Delta_{\SSS^{2}}=\frac{1}{2}\sum_{\substack{1 \leq j,k \leq 3 \\ j \neq k}}(v_j \partial_k-v_k \partial_j)^2,$$
on the unit sphere $\SSS^{2}$. In view of this link between the linearized Boltzmann and Landau operators, and in analogy with the Gelfand-Shilov smoothing result proven in~\cite{LMPX3} for the spatially homogeneous non-cutoff Boltzmann equation, we may therefore expect that the spatially homogeneous Landau equation also enjoys specific Gelfand-Shilov smoothing properties. The purpose of this note is to confirm this insight and to check that the Cauchy problem for the fluctuation around the Maxwellian equilibrium distribution associated to the spatially homogeneous Landau equation with Maxwellian molecules actually enjoys a Gelfand-Shilov regularizing effect in the class $S_{1/2}^{1/2}(\rr^d)$, implying the ultra-analyticity of both the fluctuation and its Fourier transform, for any positive time.

\section{The Landau equation}
The Landau equation written by Landau in 1936 \cite{landau} is the equation
\begin{equation}\label{e1}
\begin{cases}
\partial_tf+v\cdot\nabla_{x}f=Q_L(f,f),\\
f|_{t=0}=f_0,
\end{cases}
\end{equation}
for the density distribution of the particles $f=f(t,x,v) \geq 0$  at time $t$, having position $x \in \rr^d$ and velocity $v \in \rr^d$, with $d \geq 2$. The term $Q_L(f,f)$ is the Landau collision operator associated to the Landau bilinear operator
$$Q_L(g, f)=\nabla_v \cdot \Big(\int_{\RR^d}a(v-v_*)\big(g(t,x,v_*)(\nabla_v f)(t,x,v)-(\nabla_v g)(t,x,v_*)f(t,x,v)\big)d v_*\Big),$$
where $a=(a_{i,j})_{1 \leq i,j \leq d}$ stands for the non-negative symmetric matrix
\begin{equation}\label{landau_collision1}
a(v)=|v|^{\gamma}(|v|^2\Id -v\otimes v) \in M_d(\rr),  \quad -d<\gamma<+\infty.
\end{equation}
In this work, we study the spatially homogeneous case when the density distribution of the particles does not depend on the position variable
\begin{equation}\label{vkk1}
\begin{cases}
\partial_tf=Q_L(f,f),\\
f|_{t=0}=f_0,
\end{cases}
\end{equation}
for Maxwellian molecules, that is, when the parameter $\gamma=0$ in the assumption (\ref{landau_collision1}). At least formally, it is easily checked that the mass, the momentum and the kinetic energy are conserved quantities by this evolution equation
\begin{equation}\label{jen0}
\int_{\rr^d}f(t,v)dv=M, \quad \int_{\rr^d}f(t,v)vdv=MV, \quad \frac{1}{2}\int_{\rr^d}f(t,v)|v|^2dv=E, \quad t \geq 0,
\end{equation}
with $M > 0$, $V \in \rr^d$, $E > 0$. The Cauchy problem (\ref{vkk1}) associated to the spatially homogeneous Landau equation with Maxwellian molecules and some quantitative features of the solutions were thoroughly studied by Villani~\cite{villani3}. The propositions~4 and~6 of the work~\cite{villani3} show that, for each non-negative measurable initial density distribution $f_0$ having finite mass and finite energy
\begin{equation}\label{jen1}
f_0 \geq 0, \quad 0<\int_{\rr^d}f_0(v)dv=M<+\infty, \quad 0<\frac{1}{2}\int_{\rr^d}f_0(v)|v|^2dv=E<+\infty,
\end{equation}
the Cauchy problem (\ref{vkk1}) admits a unique global classical solution $f(t,v)$ defined for all $t \geq 0$. Furthermore, this solution is showed to be a non-negative bounded smooth function 
$$f(t) \geq 0, \quad f(t) \in L^{\infty}(\rr_v^d) \cap C^{\infty}(\rr^d_v),$$ 
for any positive time $t>0$.

In this work, we study a close-to-equilibrium framework. To that end, 
we consider the linearization of the spatially homogeneous Landau equation 
$$f=\mu_d+\sqrt{\mu_d}g,$$
around the Maxwellian equilibrium distribution 
\begin{equation}\label{maxwe1}
\mu_d(v)=(2\pi)^{-\frac{d}{2}}e^{-\frac{|v|^2}{2}}, \quad v \in \rr^d.
\end{equation}
By using that $Q_L(\mu_d,\mu_d)=0$, and setting
\begin{equation}\label{jan1}
\mathscr{L}_Lg=-\mu_d^{-1/2}Q(\mu_d,\mu_d^{1/2}g)-\mu_d^{-1/2}Q(\mu_d^{1/2}g,\mu_d),
\end{equation}
the original spatially homogeneous Landau equation (\ref{vkk1}) is reduced to the Cauchy problem for the fluctuation
\begin{equation}\label{boltz}
\begin{cases}
\partial_tg+\mathscr{L}_Lg=\mu_d^{-1/2}Q_L(\sqrt{\mu_d}g,\sqrt{\mu_d}g),\\
g|_{t=0}=g_0.
\end{cases}
\end{equation}
An explicit computation~\cite{LMPX2} (Proposition~2.1) shows that the linearized Landau operator with Maxwellian molecules acting on the Schwartz space is equal to
\begin{multline*}
\mathscr{L}_L=(d-1)\Big(\mathcal{H}-\frac{d}{2}\Big)-\Delta_{\SSS^{d-1}}+\Big[\Delta_{\SSS^{d-1}}-(d-1)\Big(\mathcal{H}-\frac{d}{2}\Big)\Big]\mathbb{P}_1\\ +\Big[-\Delta_{\SSS^{d-1}}-(d-1)\Big(\mathcal{H}-\frac{d}{2}\Big)\Big]\mathbb{P}_2,
\end{multline*}
where $\mathcal{H}=-\Delta_v+\frac{|v|^2}{4}$ is the harmonic oscillator,
$$\Delta_{\SSS^{d-1}}=\frac{1}{2}\sum_{\substack{1 \leq j,k \leq d \\ j \neq k}}(v_j \partial_k-v_k \partial_j)^2,$$
stands for the Laplace-Beltrami operator on the unit sphere $\SSS^{d-1}$ and $\mathbb{P}_{k}$ are the orthogonal projections onto the Hermite basis defined in Section~\ref{6.sec.harmo}. 
The linearized Landau operator is a non-negative operator 
$$(\mathscr{L}_Lg,g)_{L^2(\rr^d_{v})} \geq 0,$$
satisfying 
\begin{equation}\label{ker}
(\mathscr{L}_Lg,g)_{L^2(\rr^d)}=0 \Leftrightarrow g=\mathbf{P}g,
\end{equation}
where $\mathbf{P}g=(a+b \cdot v+c|v|^2)\mu_d^{1/2},$ with $a,c \in \rr$, $b \in \rr^d$, stands for the $L^2$-orthogonal projection onto the space of collisional invariants
\begin{equation}\label{coli}
\mathcal{N}=\textrm{Span}\big\{\mu_d^{1/2},v_1 \mu_d^{1/2},...,v_d\mu_d^{1/2},|v|^2\mu_d^{1/2}\big\}.
\end{equation}

By elaborating on the solutions constructed by Villani~\cite{villani3}, the purpose of this note is to study the Gelfand-Shilov regularizing properties of the Cauchy problem (\ref{boltz}) 
for the fluctuation around the Maxwellian equilibrium distribution.
For the sake of simplicity, we may assume without loss of generality that the density distribution satisfies (\ref{jen0}) with $V=0$. Furthermore, by changing the unknown function $f$ to $\tilde{f}$ as 
\begin{equation}\label{jen2.1}
f=\frac{M}{\alpha^{d}}\tilde{f}\Big(\frac{\cdot}{\alpha}\Big), \qquad \alpha=\sqrt{\frac{2E}{Md}},
\end{equation} 
we may reduce our study to the case when 
\begin{equation}\label{jen1.1}
\int_{\rr^d}f(t,v)dv=1, \quad \int_{\rr^d}f(t,v)vdv=0, \quad \int_{\rr^d}f(t,v)|v|^2dv=d, \quad t \geq 0.
\end{equation}
Let $f_0=\mu_d+\sqrt{\mu_d}g_0 \geq 0$, with $g_0 \in L^1(\rr_v^d) \cap L^2(\rr_v^d)$, be a non-negative initial density distribution having finite mass and finite energy such that
\begin{equation}\label{jen2.3}
\int_{\rr^d}f_0(v)dv=1,  \quad \int_{\rr^d}f_0(v)vdv=0, \quad \int_{\rr^d}f_0(v)|v|^2dv=d.
\end{equation}
Such an initial density distribution is rapidly decreasing with a finite temperature tail
$$\frac{1}{2} \leq \frac{1}{T(f_0)}=\sup\Big\{\beta \geq 0 : \int_{\rr^d}f_0(v)e^{\beta\frac{|v|^2}{2}}dv<+\infty\Big\},$$
since
\begin{equation}\label{j2}
\int_{\rr^d}f_0(v)e^{\frac{|v|^2}{4}}dv=\frac{1}{(2 \pi)^{\frac{d}{4}}}\int_{\rr^d}\big(\sqrt{\mu_d(v)}+g_0(v)\big)dv<+\infty,
\end{equation}
when $g_0 \in L^1(\rr_v^d)$. The analysis of the evolution of the temperature tail led in~\cite{villani3} (Section~6, p.~972-974) shows that 
$$\int_{\rr^d}f_0(v)e^{\frac{|v|^2}{4}}dv <+\infty \Rightarrow \forall t > 0, \quad \int_{\rr^d}f(t,v)e^{\frac{|v|^2}{4}}dv<+\infty.$$
This implies that the fluctuation $f=\mu_d+\sqrt{\mu_d}g \geq 0,$ around the Maxwellian equilibrium distribution defined by
\begin{equation}\label{j1}
g(t)=\mu_d^{-1/2}(f(t)-\mu_d) \in L^1(\rr_v^d) \cap C^{\infty}(\rr_v^d) \subset \mathscr{S}'(\rr_v^d), \quad t>0,
\end{equation}  
belongs to $L^1(\rr_v^d)$ and therefore remains a tempered distribution for all $t>0$. The following statement is the main result contained in this note:

\bigskip

\begin{theorem}\label{th1}
Let $f_0=\mu_d+\sqrt{\mu_d}g_0 \geq 0$, with $g_0 \in L^1(\rr_v^d) \cap L^2(\rr_v^d)$, be a non-negative measurable function having finite mass and finite energy such that
\begin{equation}\label{jen2.3}
\int_{\rr^d}f_0(v)dv=1,  \quad \int_{\rr^d}f_0(v)vdv=0, \quad \int_{\rr^d}f_0(v)|v|^2dv=d.
\end{equation}
Let $f(t)=\mu_d+\sqrt{\mu_d}g(t)$, with $g(t)  \in L^1(\rr_v^d) \cap C^{\infty}(\rr_v^d)$ when $t > 0$, be the unique global classical solution of the Cauchy problem associated to the spatially homogeneous Landau equation with Maxwellian molecules
$$\begin{cases}
\partial_tf=Q_L(f,f),\\
f|_{t=0}=f_0,
\end{cases}$$
constructed by Villani~\cite{villani3}. Then, there exists a positive constant $\delta>0$ such that 
$$\exists C>0, \forall t \geq 0, \quad \|e^{t\delta \mathcal{H}}g(t)\|_{L^2}=\Big(\sum_{k\ge 0}e^{\delta (2k+d)t}\|\mathbb P_{k}g(t)\|_{L^2}^2\Big)^{1/2} \leq Ce^{d(d-1)t}(\|g_0\|_{L^2}+1),$$
with $\mathcal{H}=-\Delta_v+\frac{|v|^2}{4},$ 
where $\|\cdot\|_{L^2}$ stands for the $L^2(\rr_v^d)$-norm and $\mathbb{P}_{k}$ are the orthogonal projections onto the Hermite basis defined in Section~\ref{6.sec.harmo}. 
In particular, this implies that the fluctuation belongs to the Gelfand-Shilov space $S_{1/2}^{1/2}(\rr^d)$ for any positive time
$$\forall t>0, \quad g(t) \in S_{1/2}^{1/2}(\rr^d).$$
\end{theorem}

\bigskip

\noindent
\textbf{Remark.}
The orthogonal projection $\mathbb{P}_{k} : \mathscr{S}'(\rr_v^d) \rightarrow \mathscr{S}(\rr_v^d)$ is well-defined on tempered distributions since the Hermite functions are Schwartz functions. 

\bigskip

This result shows that the Cauchy problem (\ref{boltz}) enjoys an ultra-analytic regularizing effect in the Gevrey class $G^{1/2}(\rr^d)$ both for the fluctuation and its Fourier transform in the velocity variable for any positive time
$$g(t),\ \widehat{g}(t) \in G^{1/2}(\rr^d), \quad t>0.$$
Let us recall that the existence, uniqueness, the Sobolev regularity and the polynomial decay of the weak solutions to the Cauchy problem (\ref{vkk1}) have been studied by Desvillettes and Villani for hard potentials~\cite{desvill} (Theorem~6), that is, when the parameter satisfies $0<\gamma \leq 1$ in the assumption (\ref{landau_collision1}). Under rather weak assumptions on the initial datum, e.g. $f_0 \in L^1_{2+\delta}$, with $\delta>0$, they prove that there exists a weak solution to the Cauchy problem such that $f \in C^{\infty}([t_0,+\infty[,\mathscr{S}(\rr_v^d))$, for all $t_0>0$, and for all $t_0 >0$, $s>0$, $m \in \nn$, 
$$\sup_{t \geq t_0}\|f(t,\cdot)\|_{H_s^m} <+\infty.$$
The Gevrey regularity $f(t,\cdot) \in G^{\sigma}$, for any $\sigma>1$, for all positive time $t>0$ of the solution to the Cauchy problem (\ref{vkk1}) with an initial datum $f_0$ with finite mass, energy and entropy satisfying
$$\forall t_0>0, \ m \geq 0, \quad \sup_{t \geq t_0}\|f(t,\cdot)\|_{H_\gamma^m} <+\infty, $$
was later established by Chen, Li and Xu for the hard potential case and the Maxwellian molecules case~\cite{CLX1}. Under the same assumptions on the solution, this result was later extended to analytic regularity~\cite{CLX2}:
$$\forall t_0>0, \exists c_0, C>0, \forall t \geq t_0, \quad \|e^{c_0(-\Delta_v)^{1/2}}f(t,\cdot)\|_{L^2} \leq C(t+1),$$   
in the hard potential case and the Maxwellian molecules case. Regarding specifically the Maxwellian molecules case $\gamma=0$, Morimoto and Xu established in the ultra-analyticity~\cite{M-X2} (Theorem~1.1), 
$$\forall \ 0<t<T, \quad f(t,\cdot) \in G^{1/2}(\rr^d),$$ 
$$\forall \ 0<T_0<T, \exists c_0>0, \forall \ 0<t\leq T_0, \quad \|e^{-c_0t \Delta_v}f(t,\cdot)\|_{L^2} \leq e^{\frac{d}{2}t}\|f_0\|_{L^2},$$
of any positive weak solution $f(t,x)>0$  to the Cauchy problem (\ref{vkk1}) satisfying $f \in L^{\infty}(]0,T[,L^2(\rr^d)\cap L_2^1(\rr^d))$, with $0<T \leq +\infty$, with an initial datum satisfying $f_0 \in L^2(\rr^d) \cap L_2^1(\rr^d)$.
The result of Theorem~\ref{th1} allows to specify further the property of ultra-analytic smoothing proven by Morimoto and Xu in the close-to-equilibrium framework~\cite{M-X2}. This result points out the specific decay of the fluctuation both in the velocity and its dual Fourier variable. As for the Boltzmann equation, the Gelfand-Shilov regularity seems relevant to describe the regularizing properties of the Landau equation in the close-to-equilibrium framework.

\section{Proof of Theorem~\ref{th1}}
The proof of Theorem~\ref{th1} is elementary and relies only on spectral arguments following the results established by Villani~\cite{villani3}.
Let $f_0=\mu_d+\sqrt{\mu_d}g_0 \geq 0$, with $g_0 \in L^1(\rr_v^d) \cap L^2(\rr_v^d)$, be a non-negative measurable function having finite mass and finite energy such that 
\begin{equation}\label{jen2}
\int_{\rr^d}f_0(v)dv=1,  \quad \int_{\rr^d}f_0(v)vdv=0, \quad \int_{\rr^d}f_0(v)|v|^2dv=d.
\end{equation}
Following~\cite{villani3} (p.~966), we may choose an orthonormal basis of $\rr^d$ diagonalizing the non-negative symmetric quadratic form
$$q(x)=\int_{\rr^d}f_0(v)(x \cdot v)^2dv=\sum_{j,k=1}^dx_j x_k \int_{\rr^d}f_0(v)v_jv_kdv \geq 0,$$ 
where $x \cdot v =\sum_{j=1}^dx_j v_j$, $x=(x_1,...,x_d)$, $v=(v_1,...,v_d)$, stands for the standard dot product in $\rr^d$.
In this orthonormal basis of $\rr_v^d$, the unique solution to the Cauchy problem associated to the spatially homogeneous Landau equation with Maxwellian molecules (\ref{vkk1}) is showed to satisfy~\cite{villani3} (Section~5), 
\begin{equation}\label{jen4}
\partial_tf=\sum_{j=1}^d(d-T_j(t))\partial_{j}^2f+(d-1)\nabla \cdot (vf)+\Delta_{\SSS^{d-1}}f,
\end{equation}
with
$$T_j(t)=\int_{\rr^d}f(t,v)v_j^2dv=1+(T_j(0)-1)e^{-4dt},$$
and
\begin{multline}\label{j5}
\int_{\rr^d}f(t,v)dv=1,  \quad \int_{\rr^d}f(t,v)v_jdv=0, \quad \int_{\rr^d}f(t,v)|v|^2dv=\sum_{j=1}^dT_j(t)=d,\\
j \neq k \Rightarrow \int_{\rr^d}f(t,v)v_jv_k=0,
\end{multline}
when $t \geq 0$.
These conditions imply that the fluctuation satisfies $g(t) \in L^1(\rr_v^d) \cap \mathcal{N}^{\perp}$, that is, 
$$\int_{\rr^d}\sqrt{\mu_d(v)}g(t,v)dv=0,  \quad \int_{\rr^d}v_j\sqrt{\mu_d(v)}g(t,v)dv=0, \quad \int_{\rr^d}|v|^2\sqrt{\mu_d(v)}g(t,v)dv=0,$$
together with
$$j \neq k \Rightarrow \int_{\rr^d}v_jv_k\sqrt{\mu_d(v)}g(t,v)=0,$$
when $t \geq 0$.
The equation (\ref{jen4}) may be rewritten for the fluctuation as
\begin{multline*}
\partial_tg=\mu_d^{-1/2}\sum_{j=1}^d(d-1-\alpha_je^{-4dt})\partial_j^2(\mu_d+\sqrt{\mu_d}g)+(d-1)\mu_d^{-1/2} \nabla \cdot (v\mu_d+v\sqrt{\mu_d}g)+\Delta_{\SSS^{d-1}}g,
\end{multline*}
with
\begin{equation}\label{kk12}
\alpha_j=\int_{\rr^d}v_j^2\sqrt{\mu_d(v)}g_0(v)dv, \quad \sum_{j=1}^d\alpha_j=0.
\end{equation}
It follows that 
\begin{multline*}
\partial_tg=-\Big[(d-1)\Big(-\Delta_v+\frac{|v|^2}{4}-\frac{d}{2}\Big)-\Delta_{\SSS^{d-1}}\Big]g-e^{-4dt}\sum_{j=1}^d\alpha_j\Big[\partial_j^2+\frac{v_j^2}{4}-v_j\partial_j-\frac{1}{2}\Big]g\\
-e^{-4dt}\sum_{j=1}^d\alpha_j(v_j^2-1)\mu_d^{1/2}.
\end{multline*}
By using that $\sum_{j=1}^d\alpha_j=0$, we notice that 
\begin{multline*}
\partial_tg=-\Big[(d-1)\Big(-\Delta_v+\frac{|v|^2}{4}-\frac{d}{2}\Big)-\Delta_{\SSS^{d-1}}\Big]g-e^{-4dt}\sum_{j=1}^d\alpha_j[(A_{+,j})^2g+v_j^2\mu_d^{1/2}],
\end{multline*}
where $A_{+,j}$ is the creation operator defined in Section~\ref{6.sec.harmo}.
We consider
$$\mathbf{S}_n=\sum_{k=0}^n\mathbb{P}_k,$$
the orthogonal projection onto the $n+1$ lowest energy levels of the harmonic oscillator, where  $\mathbb{P}_{k}$ stands for the orthogonal projection onto the Hermite basis defined in Section~\ref{6.sec.harmo}. As mentioned above, the orthogonal projection $\mathbf{S}_n$ is well-defined on tempered distributions since the Hermite functions are Schwartz functions. This gives a sense for the orthogonal projection of the fluctuation $\mathbf{S}_ng(t) \in \mathscr{S}(\rr_v^d)$ as a Schwartz function.
Then, a direct computation shows that for all $t \geq 0$, $\delta>0$, $n \geq 2$,
\begin{align*}
& \ \frac{1}{2}\partial_t(\|e^{t\delta \mathcal{H}}\mathbf{S}_ng\|_{L^2}^2)-\delta (\mathcal{H}(e^{t\delta \mathcal{H}}\mathbf{S}_ng),e^{t\delta \mathcal{H}}\mathbf{S}_ng)_{L^2}=\textrm{Re}(\partial_t\mathbf{S}_ng,e^{2\delta t\mathcal{H}}\mathbf{S}_ng)_{L^2}\\
=& \ -(d-1)(\mathcal{H}(e^{t\delta \mathcal{H}}\mathbf{S}_ng),e^{t\delta \mathcal{H}}\mathbf{S}_ng)_{L^2}
-((-\Delta_{\SSS^{d-1}})(e^{t\delta \mathcal{H}}\mathbf{S}_ng),e^{t\delta \mathcal{H}}\mathbf{S}_ng)_{L^2}\\
& \ +\frac{1}{2}d(d-1)\|e^{t\delta \mathcal{H}}\mathbf{S}_ng\|_{L^2}^2-e^{-4dt}\sum_{j=1}^d\alpha_j(e^{t\delta \mathcal{H}}(v_j^2\mu_d^{1/2}),e^{t\delta \mathcal{H}}\mathbf{S}_ng)_{L^2}\\
&\qquad -e^{-4dt}\sum_{j=1}^d\alpha_j(e^{t\delta \mathcal{H}}\mathbf{S}_n(A_{+,j})^2g,e^{t\delta \mathcal{H}}\mathbf{S}_ng)_{L^2},
\end{align*}
since the harmonic oscillator and the Laplace-Beltrami operator on $\SSS^{d-1}$ are commuting selfadjoint operators. We deduce from (\ref{kl1}), (\ref{jen6}) and (\ref{jen5}) that 
\begin{multline*}
(e^{t\delta \mathcal{H}}\mathbf{S}_n(A_{+,j})^2g,e^{t\delta \mathcal{H}}\mathbf{S}_ng)_{L^2}=e^{2\delta t}((A_{+,j})^2e^{t\delta \mathcal{H}}\mathbf{S}_{n-2}g,e^{t\delta \mathcal{H}}\mathbf{S}_ng)_{L^2}\\
=e^{2\delta t}(A_{+,j}e^{t\delta\mathcal{H}}\mathbf{S}_{n-2}g,A_{-,j}e^{t\delta \mathcal{H}}\mathbf{S}_ng)_{L^2}
\end{multline*}
and 
\begin{multline*}
e^{t\delta \mathcal{H}}(v_j^2\mu_d^{1/2})=e^{t\delta \mathcal{H}}((A_{+,j}+A_{-,j})^2\Psi_0)=e^{t\delta \mathcal{H}}(A_{+,j}^2+A_{-,j}^2+A_{+,j}A_{-,j}+A_{-,j}A_{+,j})\Psi_0\\
=(e^{2\delta t}A_{+,j}^2+e^{-2\delta t}A_{-,j}^2+A_{+,j}A_{-,j}+A_{-,j}A_{+,j})e^{\frac{d}{2}\delta t}\Psi_0=\sqrt{2}e^{(2+\frac{d}{2})\delta t}\Psi_{2e_j}+e^{\frac{d}{2}\delta t}\Psi_0.
\end{multline*}
It follows that
\begin{multline*}
 \frac{1}{2}\partial_t(\|e^{t\delta \mathcal{H}}\mathbf{S}_ng\|_{L^2}^2)+(d-1-\delta)(\mathcal{H}(e^{t\delta \mathcal{H}}\mathbf{S}_ng),e^{t\delta \mathcal{H}}\mathbf{S}_ng)_{L^2}\\ \leq \frac{1}{2}d(d-1)\|e^{t\delta \mathcal{H}}\mathbf{S}_ng\|_{L^2}^2
 +e^{-(4d-2\delta ) t}\sum_{j=1}^d|\alpha_j|\|A_{+,j}e^{t\delta \mathcal{H}}\mathbf{S}_{n-2}g\|_{L^2}\|A_{-,j}e^{t\delta \mathcal{H}}\mathbf{S}_ng\|_{L^2}\\ +e^{-(4-\frac{\delta}{2})dt}\sqrt{2e^{4\delta t}+1}\Big(\sum_{j=1}^d|\alpha_j|\Big)\|e^{t\delta \mathcal{H}}\mathbf{S}_ng\|_{L^2}.
\end{multline*}
By using that 
$$\mathcal{H}=\frac{1}{2}\sum_{j=1}^d(A_{+,j}A_{-,j}+A_{-,j}A_{+,j}),$$
we notice that 
$$\sum_{j=1}^d\|A_{+,j}u\|_{L^2}\|A_{-,j}u\|_{L^2} \leq \frac{1}{2}\sum_{j=1}^d(\|A_{+,j}u\|_{L^2}^2+\|A_{-,j}u\|_{L^2}^2)=(\mathcal{H}u,u)_{L^2}.$$
By using that
$$\|A_{+,j}e^{t\delta\mathcal{H}}\mathbf{S}_{n-2}g\|_{L^2} \leq \|A_{+,j}e^{t\delta\mathcal{H}}\mathbf{S}_{n}g\|_{L^2},$$
we obtain that 
\begin{multline*}
 \frac{1}{2}\partial_t(\|e^{t\delta \mathcal{H}}\mathbf{S}_ng\|_{L^2}^2)+(d-1-\delta)(\mathcal{H}(e^{t\delta \mathcal{H}}\mathbf{S}_ng),e^{t\delta \mathcal{H}}\mathbf{S}_ng)_{L^2}  \\ \leq  \frac{1}{2}d(d-1)\|e^{t\delta \mathcal{H}}\mathbf{S}_ng\|_{L^2}^2
+ e^{-(4d-2\delta)t}\Big(\sup_{1 \leq j \leq d}|\alpha_j|\Big)(\mathcal{H}(e^{t\delta \mathcal{H}}\mathbf{S}_ng),e^{t\delta \mathcal{H}}\mathbf{S}_ng)_{L^2}\\ 
 +e^{-(4-\frac{\delta}{2})dt}\sqrt{2e^{4\delta t}+1}\Big(\sum_{j=1}^d|\alpha_j|\Big)\|e^{t\delta \mathcal{H}}\mathbf{S}_ng\|_{L^2}.
\end{multline*}
We notice that 
$$0 < 1+\alpha_j=\int_{\rr^d}v_j^2(\mu_d+\sqrt{\mu_d}g_0)dv,$$
since the initial density distribution $f_0=\mu_d+\sqrt{\mu_d}g_0 \geq 0$ satisfies 
$\int_{\rr^d}f_0dv=1$.
On the other hand, we deduce from (\ref{kk12}) that
$$\sum_{j=1}^d(1+\alpha_j)=d.$$
This implies that $-1 < \alpha_j < d-1$, because $d \geq 2$. We may choose the positive constant $0<\delta \leq 1$ such that
$$\sup_{1 \leq j \leq d}|\alpha_j| \leq d-1-\delta.$$
It follows that 
\begin{multline*}
 \frac{1}{2}\partial_t(\|e^{t\delta \mathcal{H}}\mathbf{S}_{n}g\|_{L^2}^2)  \leq 
 \frac{1}{2}d(d-1)\|e^{t\delta\mathcal{H}}\mathbf{S}_{n}g\|_{L^2}^2
 +e^{-(4-\frac{\delta}{2})dt}\sqrt{2e^{4\delta t}+1}\Big(\sum_{j=1}^d|\alpha_j|\Big)\|e^{t\delta \mathcal{H}}\mathbf{S}_{n}g\|_{L^2}\\
\leq 
 \frac{1}{2}d(d-1)\|e^{t\delta\mathcal{H}}\mathbf{S}_{n}g\|_{L^2}^2
 +\sqrt{3}d(d-1)\|e^{t\delta \mathcal{H}}\mathbf{S}_{n}g\|_{L^2} \leq  
 d(d-1)\|e^{t\delta\mathcal{H}}\mathbf{S}_{n}g\|_{L^2}^2+\frac{9}{2}d(d-1). 
\end{multline*}
We obtain that for all $t \geq 0$, $n \geq 2$,
$$\|e^{t\delta \mathcal{H}}\mathbf{S}_{n}g(t)\|_{L^2}^2 \leq e^{2d(d-1)t}\|g_0\|_{L^2}^2+\frac{9}{2}(e^{2d(d-1)t}-1),$$
which implies that for all $t \geq 0$,
$$\|e^{t\delta \mathcal{H}}g(t)\|_{L^2}^2 \leq e^{2d(d-1)t}\|g_0\|_{L^2}^2+\frac{9}{2}(e^{2d(d-1)t}-1).$$
It follows that there exists a positive constant $C>0$ such that 
$$\|e^{t\delta \mathcal{H}}g(t)\|_{L^2(\rr_v^d)} \leq Ce^{d(d-1)t}(\|g_0\|_{L^2}+1), \quad t \geq 0,$$
and we deduce from (\ref{gel4}) that for any positive time
$$g(t) \in S_{1/2}^{1/2}(\rr^d), \quad t>0.$$
This ends the proof of Theorem~\ref{th1}.

\section{Appendix}\label{appendix}
\subsection{The harmonic oscillator}\label{6.sec.harmo}
The standard Hermite functions $(\phi_{n})_{n\geq 0}$ are defined for $x \in \rr$,
 \begin{multline}
 \phi_{n}(x)=\frac{(-1)^n}{\sqrt{2^n n!\sqrt{\pi}}} e^{\frac{x^2}{2}}\frac{d^n}{dx^n}(e^{-x^2})
 =\frac{1}{\sqrt{2^n n!\sqrt{\pi}}} \Bigl(x-\frac{d}{dx}\Bigr)^n(e^{-\frac{x^2}{2}})=\frac{ a_{+}^n \phi_{0}}{\sqrt{n!}},
\end{multline}
where $a_{+}$ is the creation operator
$$a_{+}=\frac{1}{\sqrt{2}}\Big(x-\frac{d}{dx}\Big).$$
The family $(\phi_{n})_{n\geq 0}$ is an orthonormal basis of $L^2(\R)$.
We set for $n\geq 0$, $\alpha=(\alpha_{j})_{1\le j\le d}\in\N^d$, $x\in \R$, $v\in \R^d,$
\begin{align}\label{}
\psi_n(x)&=2^{-1/4}\phi_n(2^{-1/2}x),\quad \psi_{n}=\frac{1}{\sqrt{n!}}\Bigl(\frac{x}2-\frac{d}{dx}\Bigr)^n\psi_{0},
\\
\Psi_{\alpha}(v)&=\prod_{j=1}^d\psi_{\alpha_j}(v_j),\quad \mathcal E_{k}=\text{Span}
\{\Psi_{\alpha}\}_{\alpha\in \N^d,\val \alpha=k},
\end{align}
with $\val \alpha=\alpha_{1}+\dots+\alpha_{d}$. The family $(\Psi_{\alpha})_{\alpha \in \nn^d}$ is an orthonormal basis of $L^2(\R^d)$
composed by the eigenfunctions of the $d$-dimensional harmonic oscillator
\begin{equation}\label{6.harmo}
\mathcal{H}=-\Delta_v+\frac{|v|^2}{4}=\sum_{k\ge 0}\Big(\frac d2+k\Big)\mathbb P_{k},\quad \text{Id}=\sum_{k \ge 0}\mathbb P_{k},
\end{equation}
where $\mathbb P_{k}$ is the orthogonal projection onto $\mathcal E_{k}$ whose dimension is $\binom{k+d-1}{d-1}$. The eigenvalue
$d/2$ is simple in all dimensions and $\mathcal E_{0}$ is generated by the function
\begin{equation}\label{kl1}
\Psi_{0}(v)=\frac{1}{(2\pi)^{\frac{d}{4}}}e^{-\frac{\val v^2}{4}}=\mu_d^{1/2}(v),
\end{equation}
where $\mu_d$ is the Maxwellian distribution defined in (\ref{maxwe}).
Setting
$$A_{\pm,j}=\frac{ v_{j}}2\mp\frac{\p}{\p v_{j}}, \quad 1 \leq j \leq d,$$
we have
$$\Psi_{\alpha}=\frac{1}{\sqrt{\alpha_{1}!... \alpha_{d}!}}A_{+,1}^{\alpha_{1}}... A_{+,d}^{\alpha_{d}} \Psi_{0}, \quad \alpha=(\alpha_1,...,\alpha_d)\in \N^d,$$
\begin{equation}\label{jen6}
A_{+,j}\Psi_{\alpha}=\sqrt{\alpha_{j}+1}\Psi_{\alpha+e_{j}},\quad A_{-,j}\Psi_{\alpha}=\sqrt{\alpha_{j}}\Psi_{\alpha-e_{j}} \ (=0 \textrm{ if } \alpha_j=0),
\end{equation}
where $(e_1,...,e_d)$ stands for the canonical basis of $\rr^d$. In particular, we readily notice that for all $t \geq 0$, $\delta>0$,
\begin{equation}\label{jen5}
e^{t\delta \mathcal{H}}A_{+,j}=e^{\delta t}A_{+,j}e^{t\delta \mathcal{H}}, \quad e^{t\delta \mathcal{H}}A_{-,j}=e^{-\delta t}A_{-,j}e^{t\delta \mathcal{H}}.
\end{equation}

\subsection{Gelfand-Shilov regularity}\label{regularity}We refer the reader to the works~\cite{gelfand,rodino1,rodino,toft} and the references herein for extensive expositions of the Gelfand-Shilov regularity.  The Gelfand-Shilov spaces $S_{\nu}^{\mu}(\rr^d)$, with $\mu,\nu>0$, $\mu+\nu\geq 1$, are defined as the spaces of smooth functions $f \in C^{\infty}(\rr^d)$ satisfying to the estimates
$$\exists A,C>0, \quad |\partial_v^{\alpha}f(v)| \leq C A^{|\alpha|}(\alpha !)^{\mu}e^{-\frac{1}{A}|v|^{1/\nu}}, \quad v \in \rr^d, \ \alpha \in \nn^d,$$
or, equivalently
$$\exists A,C>0, \quad \sup_{v \in \rr^d}|v^{\beta}\partial_v^{\alpha}f(v)| \leq C A^{|\alpha|+|\beta|}(\alpha !)^{\mu}(\beta !)^{\nu}, \quad \alpha, \beta \in \nn^d.$$
These Gelfand-Shilov spaces  $S_{\nu}^{\mu}(\rr^d)$ may also be characterized as the spaces of Schwartz functions $f \in \mathscr{S}(\rr^d)$ satisfying to the estimates
$$\exists C>0, \eps>0, \quad |f(v)| \leq C e^{-\eps|v|^{1/\nu}}, \quad v \in \rr^d, \qquad |\widehat{f}(\xi)| \leq C e^{-\eps|\xi|^{1/\mu}}, \quad \xi \in \rr^d.$$
In particular, we notice that Hermite functions belong to the symmetric Gelfand-Shilov space  $S_{1/2}^{1/2}(\rr^d)$. More generally, the symmetric Gelfand-Shilov spaces $S_{\mu}^{\mu}(\rr^d)$, with $\mu \geq 1/2$, can be nicely characterized through the decomposition into the Hermite basis $(\Psi_{\alpha})_{\alpha \in \nn^d}$, see e.g. \cite{toft} (Proposition~1.2),
\begin{multline}\label{gel4}
f \in S_{\mu}^{\mu}(\rr^d) \Leftrightarrow f \in L^2(\rr^d), \ \exists t_0>0, \ \big\|\big((f,\Psi_{\alpha})_{L^2}\exp({t_0|\alpha|^{\frac{1}{2\mu}})}\big)_{\alpha \in \nn^d}\big\|_{l^2(\nn^d)}<+\infty\\
\Leftrightarrow f \in L^2(\rr^d), \ \exists t_0>0, \ \|e^{t_0\mathcal{H}^{1/2\mu}}f\|_{L^2}<+\infty,
\end{multline}
where $(\Psi_{\alpha})_{\alpha \in \nn^d}$ stands for the Hermite basis defined in Section~\ref{6.sec.harmo}, and where 
$$\mathcal{H}=-\Delta_v+\frac{|v|^2}{4},$$ is the $d$-dimensional harmonic oscillator.
The Cauchy problem defined by the evolution equation associated to the harmonic oscillator
\begin{equation}\label{gel1}
\begin{cases}
\partial_tf+\mathcal{H}f=0,\\
f|_{t=0}=f_0 \in L^2(\rr^d),
\end{cases}
\end{equation}
enjoys nice regularizing properties. The smoothing effect for the solutions to this Cauchy problem is naturally described in term of the Gelfand-Shilov regularity. The characterization (\ref{gel4}) proves that there is a regularizing effect for the solutions to the Cauchy problem (\ref{gel1}) in the symmetric Gelfand-Shilov space  $S_{1/2}^{1/2}(\rr^d)$ for any positive time, whereas the smoothing effect for the solutions to the Cauchy problem defined by the evolution equation associated to the fractional harmonic oscillator
\begin{equation}\label{gel2}
\begin{cases}
\partial_tf+\mathcal{H}^sf=0,\\
f|_{t=0}=f_0 \in L^2(\rr^d),
\end{cases}
\end{equation}
with $0<s<1$, occurs for any positive time in the symmetric Gelfand-Shilov space $S_{1/2s}^{1/2s}(\rr^d)$.

\bigskip

\vs\noindent
{\bf Acknowledgements.}
The research of the first author was supported by the Grant-in-Aid for Scientific Research No. 22540187, Japan Society for the Promotion of Science. The research of the second author was supported by the CNRS chair of excellence at Cergy-Pontoise University. The research of the last author was supported partially by ``The Fundamental Research Funds for Central Universities'' and the National Science Foundation of China No. 11171261.

\end{document}